\documentclass[12pt,reqno]{amsart}

\usepackage[margin=1in]{geometry}
\usepackage[T1]{fontenc}
\usepackage[utf8]{inputenc}
\usepackage{amsmath,amsfonts,amssymb,amsthm,amscd,mathtools}
\usepackage{hhline}
\usepackage{float}
\usepackage{graphicx}
\usepackage[titletoc,title]{appendix}
\usepackage[dvipsnames]{xcolor}
\usepackage{caption}
\usepackage{listings}
\usepackage{hyperref}

\usepackage[table,dvipsnames]{xcolor}
\usepackage{colortbl}

\numberwithin{equation}{section}

\newtheorem{thm}{Theorem}[section]

\newtheorem{cor}[thm]{Corollary}

\newtheorem{lem}[thm]{Lemma}
\newtheorem{prop}[thm]{Proposition}

\newtheorem{defi}[thm]{Definition}

\definecolor{deepblue}{rgb}{0,0,0.5}
\definecolor{deepred}{rgb}{0.6,0,0}
\definecolor{deepgreen}{rgb}{0,0.5,0}
\definecolor{ao}{rgb}{0.0,0.5,0.0}

\newcommand{\bburl}[1]{\textcolor{blue}{\url{#1}}}
\newcommand{\seqnum}[1]{\href{https://oeis.org/#1}{\rm\underline{#1}}}

\begin{document}

\title{Schreier Sets And the Product of Fibonacci and Catalan Numbers}

\author{Kevin Beanland}
\address{Department of Mathematics, Washington and Lee University, Lexington, VA 24450, USA}
\email{beanlandk@wlu.edu}

\author{H\`ung Vi\d{\^e}t Chu}
\address{Department of Mathematics, Washington and Lee University, Lexington, VA 24450, USA}
\email{hchu@wlu.edu}

\author{Carrie E. Finch-Smith}
\address{Department of Mathematics, Washington and Lee University, Lexington, VA 24450, USA}
\email{FinchC@wlu.edu}


\maketitle

\begin{abstract}
First, we introduce strictly super-Schreier sets, which provide a new combinatorial interpretation of products of Fibonacci and Catalan numbers. Second, we show that the numbers of sets in higher-order Schreier families do not satisfy any linear recurrence relation, thus answering a question in the existing literature. 
\end{abstract}

\begin{center}
\small{\textit{In memory of Florian Luca (1969--2026).}}
\end{center}

\tableofcontents

\section{Introduction}
A nonempty finite subset $F$ of $\mathbb{N}$ is \textit{Schreier} if $\min F\ge |F|$. The collection of Schreier sets, including the empty set, is denoted by $\mathcal{S}$ or $\mathcal{S}_1$. Schreier sets have been shown to be connected to many interesting linear recurrence sequences. For example, counting Schreier sets gives 
the Fibonacci sequence \cite{bird2012shcreierfib}, partial sums of the Fibonacci sequence \cite{C1}, the $s$-step Fibonacci sequence \cite{CIMSZ}. Here the Fibonacci sequence $(F_n)_{n=1}^\infty$ is given by $F_1 = F_2 = 1$ and $F_n = F_{n-1} + F_{n-2}$ for all $n\ge 3$. 

More recently, \cite{CKKZ} introduced \textit{$k$-super-Schreier} sets and observed a connection with the Catalan numbers. Given $i\in\mathbb{N}$ and a set $F$ with $|F|\ge i$, let $\min_i F$ denote the $i$th smallest element of $F$. 
\begin{defi}\normalfont\cite[Definition 1.3]{CKKZ} A finite set $F\subset\mathbb{N}$ is said to be \textit{$k$-super-Schreier} if $|F|\ge k$, and 
$$\min_i F\ \ge\ |F| + 2(i-1), \mbox{ for all integers }i\in [1,k].$$
\end{defi}
For each $k\in\mathbb{N}$, the $k$th Catalan number $C_k$ is
the number of lattice paths from $(0, 0)$ to $(k, k)$ which never go above the diagonal line $y=x$ \cite[Example 2.28]{KellerTrotter}, and   
$$C_k\ =\ \frac{1}{k+1}\binom{2k}{k}.$$
\begin{thm}(cf.\ \cite[Theorem 1.4]{CKKZ})
    For $k\in\mathbb{N}$, we have
    $$|\{F \subset \mathbb{N}\,:\, F \mbox{ is }k\mbox{-super-Schreier and }\max F  = 3k-2\}|\ =\ C_k.$$
\end{thm}

Catalan numbers and linear recurrences feature prominently in the work of Florian Luca, who made extensive contributions to the arithmetic theory of integer sequences. For instance, Luca and Pappalardi investigated the occurrence of binomial coefficients and Catalan numbers in binary recurrence sequences \cite{LP}, while Laishram, Luca, and Sias studied members of Lucas sequences that can be expressed as products of Catalan numbers \cite{LLS}.
More recently, Hajdu, Luca, Tengely, and Ulas examined the multiplicative properties of Catalan numbers, particularly
products that are perfect squares \cite{HLTU}.

We approach Catalan numbers and linear recurrences from an enumerative perspective. We investigate counting sequences
associated with Schreier families, their connections with Catalan numbers, and the existence or nonexistence of linear recurrence relations with constant coefficients.

The first result gives a combinatorial interpretation of the product $F_jC_k$ for every $(j,k)\in\mathbb{N}^2$.
To do so, we introduce \textit{strictly $k$-super Schreier sets}.
\begin{defi}\normalfont
    A set $F\subset\mathbb{N}$ is \textit{strictly $k$-super Schreier} if $|F|\ge k+1$, and 
    $$\min_i F\ \ge\ |F| + 2(i-1),\mbox{ for all }i\le k \quad \mbox{ and }\quad \min_{k+1} F \ <\ |F| + 2k.$$
\end{defi}

Define
$$\mathcal{R}_{k,n}\ :=\ \{F\subset \{1, 2, \ldots, n\}\,:\, F\mbox{ is strictly }k\mbox{-super Schreier}\}.$$
Table \ref{Data_oneint} reports $|\mathcal{R}_{k,n}|$ computed by a ChatGPT-written Python program \cite{AI}.
\begin{table}[H]
\centering
\begin{tabular}{ |c| c| c| c| c| c| c| c| c| c| c| c| c| c| c| c| c| c| c| c| c|c|c|c|}
\hline
$n$ &$1$& $2$ & $3$ & $4$ & $5$ & $6$ & $7$ & $8$ & $9$ & $10$ & $11$ & $12$ & $13$ & $14$ & $15$ & $16$ & $17$ & $18$\\
\hline
$|\mathcal{R}_{1,n}|$ & $0$ & $0$ & $1$ & $1$ & $2$ & $3$ & $5$ & $8$ & $13$ & $21$ & $34$ & $55$ & $89$ & $144$ & $233$ & $377$ & $610$ & $987$\\
\hline
$|\mathcal{R}_{2,n}|$ & $0$ & $0$ & $0$ & $0$ & $0$ & $2$ & $2$ & $4$ & $6$ & $10$ & $16$ & $26$ & $42$ & $68$ & $110$ & $178$ & $288$ & $466$ \\
\hline
$|\mathcal{R}_{3,n}|$ & $0$ & $0$ & $0$ & $0$ & $0$ & $0$ & $0$ & $0$ & $5$ & $5$ & $10$ & $15$ & $25$ & $40$ & $65$ & $105$ & $170$ & $275$\\
\hline
$|\mathcal{R}_{4,n}|$ & $0$ & $0$ & $0$ & $0$ & $0$ & $0$ & $0$ & $0$ & $0$ & $0$ & $0$ & $14$ & $14$ & $28$ & $42$ & $70$ & $112$ & $182$ \\
\hline
$|\mathcal{R}_{5,n}|$ & $0$ & $0$ & $0$ & $0$ & $0$ & $0$ & $0$ & $0$ & $0$ & $0$ & $0$ & $0$ & $0$ & $0$ & $42$ & $42$ & $84$ & $126$ \\
\hline

\end{tabular}
\caption{The first $18$ values of $(|\mathcal{R}_{k, n}|)_{n=1}^\infty$ with $1\le k\le 5$. Note that $(|\mathcal{R}_{1,n}|)$ is \seqnum{A000045}; $(|\mathcal{R}_{2,n}|)$ is \seqnum{A055389}; $(|\mathcal{R}_{3,n}|)$ is \seqnum{A022088}; $(|\mathcal{R}_{4,n}|)$ is \seqnum{A022348}.}
\label{Data_oneint}
\end{table}

\begin{thm}\label{M1}
    For $k\in\mathbb{N}$, we have
    $$|\mathcal{R}_{k, n}|\ =\ \begin{cases}0, &\mbox{ if } n \le 3k-1;\\ C_k F_{n-3k+1}, &\mbox{ if }n \ge 3k.\end{cases}$$
\end{thm}

The second result answers \cite[Problem 6.1]{BGHH} in the negative. Recall that the collection $k\mathcal{S}_1$ consists of all subsets of $\mathbb{N}$ that can be written as a union of at most $k$ many Schreier sets \cite{BGHH}. For each fixed $k$,  \cite[Theorem 1.1]{BGHH} establishes linear recurrences for the count of sets in $k\mathcal{S}_1$ whose maximum is $n$. The authors then asked whether there are linear recurrences for sets in $k\mathcal{S}_2$ whose maximum is $n$, where $\mathcal{S}_2$ is the convolution of $\mathcal{S}$ with itself. Precisely,
a set $F$ is in $\mathcal{S}_2$ if there exist successive nonempty sets $E_1 < E_2 < \cdots < E_\ell, E_i\in\mathcal{S}_1$, such that 
$$\bigcup_{i=1}^\ell E_i \ =\ F\quad \mbox{ and }\quad \min E_1\ \ge\ \ell.$$
In general, one can define each subsequent collection $\mathcal{S}_k$, $k\ge 1$, inductively \cite{AT}:  
$\mathcal{S}_1$ consists of Schreier sets, and once $\mathcal{S}_k$ has been defined, the collection $\mathcal{S}_{k+1}$ consists of unions of 
successive nonempty sets $E_1 < E_2 < \cdots < E_\ell, E_i\in\mathcal{S}_{k}$, with
$\min E_1\ge \ell$. 

Once $\mathcal{S}_k$ is defined, the collection $j\mathcal{S}_k$ consists of all subsets of $\mathbb{N}$ that can be written as a union of at most $j$ many sets in $\mathcal{S}_k$.

The following theorem gives a condition for when a sequence has no linear recurrence.

\begin{thm}\label{m2}
Let $(a_n)_{n\ge 1}$ be a sequence of nonnegative real numbers. Suppose
there exists $\rho>0$ such that
$$\lim_{n\to\infty} a_n^{1/n}\ =\ \rho$$
and
$$\lim_{n\to\infty}\frac{a_n}{\rho^n}\ =\ 0.$$
Then $(a_n)$ does not satisfy any eventual linear recurrence with constant
coefficients.
\end{thm}

\begin{cor}\label{mc2}
Let $k\ge 2$, and define
$$a_{k,n}\ =\ 
\left|\left\{F\in\mathcal S_k\,:\,\max F=n\right\}\right|,\mbox{ for each } n\ge 1.$$
Then $(a_{k,n})_{n\ge 1}$ does not satisfy any eventual linear recurrence
with constant coefficients. 
\end{cor}

The following corollary answers \cite[Problem 6.1]{BGHH} in the negative.

\begin{cor}\label{mc3}
Let $j\in\mathbb{N}$ and $k\ge 2$. Define 
$$b_{j, k,n}\ =\ 
\left|\left\{F\in j\mathcal S_k\, :\,\max F=n\right\}\right|,\mbox{ for each } n\ge 1.$$
Then $(b_{j, k,n})_{n\ge 1}$ does not satisfy any eventual linear recurrence
with constant coefficients. 
\end{cor}

%%%%%%%%%%%%%%%%%%%%%%%%%%%%%%%%%%%%%%%%%%%%%%%%%%%%%%%%%%%%%%%%%%%%%%%%%%%%%%%%%%%%%%%%%%%%%%%%%%%%%%%%%%
%%%%%%%%%%%%%%%%%%%%%%%%%%%%%%%%%%%%%%%%%%%%%%%%%%%%%%%%%%%%%%%%%%%%%%%%%%%%%%%%%%%%%%%%%%%%%%%%%%%%%%%%%%
%%%%%%%%%%%%%%%%%%%%%%%%%%%%%%%%%%%%%%%%%%%%%%%%%%%%%%%%%%%%%%%%%%%%%%%%%%%%%%%%%%%%%%%%%%%%%%%%%%%%%%%%%%
\section*{Declaration of generative AI use}
The authors conceived the present project and formulated both main results in this paper independently. During preparation of the manuscript, OpenAI’s ChatGPT \cite{AI} was used (i) to generate draft proofs of Theorem \ref{m2} and related propositions and lemmas in Section \ref{sect_nolinear}; (ii) to assist in checking mathematical arguments, references, notation, and exposition; and (iii) to generate the Python program used to produce the numerical data in Table 1. The proof of Theorem \ref{M1}, including the related propositions and lemmas, were developed entirely by the authors without AI assistance. The authors independently verified every AI-assisted argument, corrected and revised the generated material as necessary, and take full responsibility for the content of the paper. 
%%%%%%%%%%%%%%%%%%%%%%%%%%%%%%%%%%%%%%%%%%%%%%%%%%%%%%%%%%%%%%%%%%%%%%%%%%%%%%%%%%%%%%%%%%%%%%%%%%%%%%%%%%%
%%%%%%%%%%%%%%%%%%%%%%%%%%%%%%%%%%%%%%%%%%%%%%%%%%%%%%%%%%%%%%%%%%%%%%%%%%%%%%%%%%%%%%%%%%%%%%%%%%%%%%%%%%%
%%%%%%%%%%%%%%%%%%%%%%%%%%%%%%%%%%%%%%%%%%%%%%%%%%%%%%%%%%%%%%%%%%%%%%%%%%%%%%%%%%%%%%%%%%%%%%%%%%%%%%%%%%%
%%%%%%%%%%%%%%%%%%%%%%%%%%%%%%%%%%%%%%%%%%%%%%%%%%%%%%%%%%%%%%%%%%%%%%%%%%%%%%%%%%%%%%%%%%%%%%%%%%%%%%%%%%%
\section{Product of Fibonacci and Catalan numbers}

To establish the product of Fibonacci and Catalan numbers in Theorem \ref{M1}, it suffices to show that $|\mathcal{R}_{k,3k}| = |\mathcal{R}_{k,3k+1}| = C_k$ and that the sequence $(|\mathcal{R}_{k,n}|)_{n=3k}^\infty$ follows the Fibonacci recurrence.

\begin{prop}\label{p1}
For $n\le 3k-1$,  we have the collection $\mathcal{R}_{k,n}$ is empty. 
\end{prop}
\begin{proof}
    Suppose that there is an $F\in \mathcal{R}_{k, n}$. We have
    $$\min_{k} F\ \ge\ |F| + 2(k-1)\ \ge\ (k+1) + 2(k-1)\ =\ 3k-1.$$
    Hence, 
    $$\min_{k+1} F \ \ge\ \min_{k} F + 1 \ \ge\ 3k,$$
    which cannot happen because $F$ is a subset of $\{1, 2, \ldots, n\}$ and $n\le 3k-1$. 
\end{proof}

\begin{prop}\label{p2}
For each $k\in\mathbb{N}$, we have
$$|\mathcal{R}_{k, 3k}| \ =\ C_k.$$
\end{prop}

\begin{proof}
    Let $F\in \mathcal{R}_{k, 3k}$. We have
    $$3k\ \ge\ \max F\ \ge\ \min_{k} F + 1\ \ge\ |F| + 2(k-1) + 1\ \ge\ (k+1) + 2k-1\ =\ 3k.$$
    Hence, 
    $$\max F \ =\ 3k, \quad |F| \ =\ k+1, \quad\mbox{ and }\quad \min_k F\ =\  3k-1.$$
    As a result, the set $F$ can be written as
    $$\{a_1 \ <\ a_2 \ <\ \cdots \ <\ a_k = 3k-1\ <\ a_{k+1} = 3k\}.$$
    
    Note that by above,
    $$\min_k F\ =\ 3k-1\ =\ |F| + 2(k-1) \quad\mbox{ and }\quad \min_{k+1} F\ =\ 3k \ <\ |F| + 2k.$$
    Hence, for $F$ to be strictly $k$-super-Schreier, we need 
    \begin{equation}\label{e5}a_1 \ \ge\ k+1, \quad a_2\ \ge\ k+3, \quad\ldots, \quad a_{k-1}\ \ge\ 3k-3.\end{equation}
    Let 
    $$x_1 \ =\ a_1 - k, \quad x_2\ =\ a_2 - a_1 - 1, \quad \ldots, \quad x_k\ =\ (3k-1)-a_{k-1}-1.$$
    Then $x_1,\dots,x_k$ are nonnegative integers satisfying
    $$x_1 + x_2 + \cdots + x_k\ =\ k,$$
    and due to \eqref{e5},
    \begin{equation}\label{e6}
        x_1 + x_2 + \cdots + x_i\ \ge\ i, \quad \mbox{ for all }1\le i\le k-1.
    \end{equation}
    Therefore, each solution $(x_1, x_2, \ldots, x_k)$ corresponds to a lattice path on the plane from $(0,0)$ to $(k, k)$ that consists of northward and eastward steps, and $x_i$ is the number of eastward steps before the $i$th northward step. Due to \eqref{e6}, such a path must not rise above the diagonal. Hence, $|\mathcal{R}_{k,3k}| = C_k$.
\end{proof}

\begin{prop}\label{pp3}
For each $k\in\mathbb{N}$, we have
$$\mathcal{R}_{k,3k}\ =\ \mathcal{R}_{k, 3k+1},$$
and thus,
$$|\mathcal{R}_{k, 3k+1}| \ =\ C_k.$$
\end{prop}

\begin{proof}
    By definition, $\mathcal{R}_{k, 3k}\subset \mathcal{R}_{k, 3k+1}$. We show that $\mathcal{R}_{k, 3k+1}\backslash \mathcal{R}_{k, 3k} = \emptyset$. Suppose that $F\in \mathcal{R}_{k, 3k+1}\backslash \mathcal{R}_{k, 3k}$. Then 
    $\max F  = 3k + 1$.

    If $|F| \ge k+2$, then
    $$3k+1 \ =\ \max F\ \ge\ \min_k F + 2\ \ge\ |F| + 2(k-1) + 2\ \ge\ 3k+2,$$
    a contradiction. Hence, $|F| = k+1$, so $\min_{k+1} F = \max F = 3k+1$. Thus, we have
    $$3k + 1\ =\ \min_{k+1} F \ <\ |F| + 2k\ =\ (k+1) + 2k\ =\ 3k+1,$$
    a contradiction. 

Therefore, $\mathcal{R}_{k, 3k+1}\backslash \mathcal{R}_{k, 3k} = \emptyset$. The second statement follows from Proposition \ref{p2}.
\end{proof}

\begin{thm}\label{ms}
For $k\in\mathbb{N}$, we have
$$|\mathcal{R}_{k, n}| + |\mathcal{R}_{k,n+1}| \ =\ |\mathcal{R}_{k, n+2}|, \quad \mbox{ whenever }n\ge 3k.$$
\end{thm}

\begin{proof}
    Define the map $f: \mathcal{R}_{k, n}\rightarrow\mathcal{R}_{k, n+2}\backslash \mathcal{R}_{k,n+1}$ as
    $$F\ \mapsto\ (F + 1)\cup \{n+2\}.$$
    It suffices to show that $f$ is a bijection. 

    First, the map $f$ is well-defined. Pick $F\in \mathcal{R}_{k,n}$. For each $i\le k$, we have
    $$\min_i f(F)\ =\ \min_i F + 1\ \ge\ |F| + 2(i-1) + 1\ =\ |f(F)| + 2(i-1).$$
    Furthermore,
    $$\min_{k+1} f(F)\ =\ \min_{k+1} F + 1\ <\ |F| + 2k + 1\ =\ |f(F)| + 2k.$$
    Hence, $f(F)$ is strictly $k$-super Schreier. Since $\max f(F) = n+2$, we have $f(F)\in \mathcal{R}_{k, n+2}\backslash \mathcal{R}_{k,n+1}$.

    The map $f$ is injective by definition. We show that $f$ is surjective. Let $G\in \mathcal{R}_{k, n+2}\backslash \mathcal{R}_{k,n+1}$. Let $E = (G\backslash \{n+2\})-1$. Then $\max E\le n$. 
    Suppose, for a contradiction, that $|G| = k+1$ and thus, $\min_{k+1} G = n+2$. Then 
    $$3k+2\ \le\ n+2\ =\ \min_{k+1} G\ <\ |G| + 2k\ =\ 3k+1,$$
    a contradiction. Hence, $|G|\ge k+2$, and so, $|E|\ge k+1$. 

We have, for $i\le k$,
    $$\min_i E\ =\ \min_i G - 1\ \ge\ |G| + 2(i-1) - 1\ =\ |E| + 2(i-1),$$
    and 
    $$\min_{k+1} E\ =\ \min_{k+1} G - 1\ <\ |G| + 2k - 1\ =\ |E| + 2k.$$
Therefore, $E$ is strictly $k$-super Schreier, so $E\in \mathcal{R}_{k,n}$. 

Since $f(E) = G$, the map $f$ is surjective. 
\end{proof}

\begin{proof}[Proof of Theorem \ref{M1}]
    Theorem \ref{M1} follows from Propositions \ref{p1}, \ref{p2}, and \ref{pp3} and Theorem \ref{ms}.
\end{proof}

\section{No linear recurrence for the union of $\mathcal{S}_k$-sets for $k\ge 2$}\label{sect_nolinear}

We start with a lemma to be used in the proof of Theorem \ref{m2}, stating that the generating function of a linear recurrence sequence must be rational. 

\begin{lem}\label{l1}
Let $(a_n)_{n\ge 1}$ be a sequence satisfying an eventual linear recurrence
with constant coefficients. Then its generating function
$G(z)=\sum_{n=1}^{\infty} a_n z^n$
is a rational function.
\end{lem}

\begin{proof}
Suppose that there exist constants $q_1,\ldots,q_d$
and $N\in\mathbb{N}$ such that
$$a_{n+d}\ =\ q_1a_{n+d-1}+\cdots+q_da_n, \mbox{ for every }n\ge N.$$
Set
$$Q(z)\ =\ 1-q_1z-\cdots-q_dz^d.$$
Then all sufficiently large coefficients of $Q(z)G(z)$ are zero, by the
recurrence relation. Hence,
$Q(z)G(z) = P(z)$ for some polynomial $P(z)$. Therefore,
$G(z) = P(z)/Q(z)$ and thus, $G$ is rational.
\end{proof}

\begin{proof}[Proof of Theorem \ref{m2}]
Suppose, toward a contradiction, that $(a_n)$ satisfies an eventual linear
recurrence with constant coefficients. By Lemma \ref{l1}, the generating
function $G(z) = \sum_{n=1}^{\infty}a_nz^n$
is rational.

Since
$\lim_{n\to\infty}a_n^{1/n} = \rho$,
the Cauchy-Hadamard formula implies that the radius of convergence of $G$
is $1/\rho$. Since $G$ is rational, it therefore has a pole $\zeta$ with 
$|\zeta|=\rho^{-1}$. Set
$$b_n\ :=\ \frac{a_n}{\rho^n}\ =\ a_n |\zeta^n|.$$
By assumption, $b_n\to 0$. We claim that
\begin{equation}\label{e21}(1-t)|G(t\zeta)|\ \rightarrow\ 0, \quad\mbox{ as }t\rightarrow 1^-.\end{equation}
Indeed, since $|\zeta|=1/\rho$,
$$|a_n(t\zeta)^n|\ =\ b_nt^n.$$
Thus,
$$(1-t)|G(t\zeta)|\ \le\ (1-t)\sum_{n=1}^{\infty}b_nt^n.$$
Given $\varepsilon>0$, choose $N$ such that $b_n\le\varepsilon$ for all
$n\ge N$. Then for $0 < t < 1$, 
\begin{align*}
(1-t)|G(t\zeta)|
&\ \le\ (1-t)\sum_{n=1}^{N-1}b_nt^n + (1-t)\sum_{n=N}^{\infty}b_nt^n\\
&\ \le\  (1-t)\sum_{n=1}^{N-1}b_n + \varepsilon(1-t)\sum_{n=N}^{\infty}t^n\\
&\ =\  (1-t)\sum_{n=1}^{N-1}b_n+\varepsilon t^N.
\end{align*}
Letting $t\to 1^-$ gives
$$\limsup_{t\to 1^-}(1-t)|G(t\zeta)|\ \le\ \varepsilon.$$
Since $\varepsilon>0$ is arbitrary,
$$(1-t)|G(t\zeta)|\ \rightarrow\ 0, \quad\mbox{ as }t\rightarrow 1^-.$$

On the other hand, suppose that the pole of $G$ at $\zeta$ has order
$d\ge 1$. Then, in a neighborhood of $\zeta$,
$$G(z)\ =\ \frac{h(z)}{(z-\zeta)^d},$$
where $h$ is analytic and $h(\zeta)\neq 0$. We have
$$(1-t)|G(t\zeta)|\ =\ \frac{|h(t\zeta)|}{|\zeta|^d}(1-t)^{1-d}, \quad \mbox{ for }0 < t < 1.$$
If $d=1$, this tends to
$$\frac{|h(\zeta)|}{|\zeta|}\ >\ 0,\quad\mbox{ as } t\rightarrow 1^-,$$
while if $d>1$, it tends to infinity. Both cases contradict \eqref{e21}. 
Therefore $(a_n)$ does not satisfy any eventual linear recurrence with
constant coefficients.
\end{proof}

We now show that certain sequences arising from counting sets in Schreier families of higher order satisfy the conditions in Theorem \ref{m2} with $\rho = 2$. To that end, we record several useful properties of the Schreier families $(\mathcal{S}_k)_{k=1}^\infty$:
\begin{enumerate}
    \item[(P1)]  By construction, $\mathcal{S}_1 \subset \mathcal{S}_2 \subset \mathcal{S}_3 \subset \cdots$, and sets in $\mathcal{S}_k$ are finite. 
    \item[(P2)] Each $\mathcal{S}_k$ is hereditary,  meaning that if $F$ is in $\mathcal{S}_k$, then all subsets of $F$ are also in $\mathcal{S}_k$ \cite[page 16]{AT}.
    \item[(P3)] \label{p3} We identify subsets of $\mathbb N$ with their characteristic functions in $\{0,1\}^{\mathbb N}$. Then each $\mathcal{S}_k$ is compact in $\{0,1\}^{\mathbb N}$ equipped with the product topology \cite[page 16]{AT}. Since $\{0,1\}^{\mathbb N}$ is a compact Hausdorff space, a subset is compact if and only if it is closed. Thus $\mathcal S_k$ is closed. In particular, if
$$F_j\ \in\ \mathcal S_k \quad\text{ and }\quad F_j\ \rightarrow\ F \mbox{ in }\{0,1\}^{\mathbb N},$$
then $F \in \mathcal S_k$.
\end{enumerate}

\begin{prop}\label{pbound}
Let $k\ge 2$, and define
$$a_{k,n}\ =\ 
\left|\left\{F\in\mathcal S_k:\max F=n\right\}\right|,\mbox{ for each } n\ge 1.$$
Then 
$$\lim_{n\rightarrow\infty} a_{k,n}^{1/n} \ =\ 2\quad\mbox{ and }\quad\lim_{n\to\infty}\frac{a_{k,n}}{2^n}\ =\ 0.$$
\end{prop}

\begin{proof}
We first show that $\lim_{n\to\infty}a_{k,n}^{1/n} = 2$.
For $m\ge 1$ and $j\in\{0,\ldots,m-1\}$, set
$$F_{m,j}\ :=\ [2^jm,2^{j+1}m-1].$$
Since
$$|F_{m,j}|\ =\ 2^jm\ =\ \min F_{m,j},$$
we have $F_{m,j}\in\mathcal S_1$, and hence,
$F_{m,j}\in\mathcal S_{k-1}$ due to Property (P1). Moreover,
$$F_{m,0}\ <\ F_{m,1}\ <\ \cdots\ <\ F_{m,m-1}$$
and the number of blocks is $m=\min F_{m,0}$.
Therefore,
$$[m,2^mm-1]\ =\ \bigcup_{j=0}^{m-1}F_{m,j}\ \in\ \mathcal S_k.$$

Choose $m = \left\lceil\log_2(n+1)\right\rceil$. Then
$$m\ \le\ n\ \le\ 2^m m-1.$$
Since $\mathcal S_k$ is hereditary (Property (P2)), every subset of $[m,n]$ containing $n$
belongs to $\mathcal S_k$. Consequently,
\begin{equation}\label{e22}2^{n-m}\ \le\ a_{k,n}\ \le\ 2^n.\end{equation}
Since $m=o(n)$, we obtain $\lim_{n\to\infty}a_{k,n}^{1/n} = 2$.

We now prove that
$$\lim_{n\to\infty}\frac{a_{k,n}}{2^n}\ =\ 0.$$
Choose a random subset $A\subset\mathbb N$ by including each positive
integer independently with probability $1/2$. For each $n\ge 1$, define
$$E_n\ =\ \{A\,:\,A\cap[1,n]\in\mathcal S_k\}.$$
Since $\mathcal S_k$ is hereditary, the events $(E_n)_{n=1}^\infty$ are decreasing. Indeed,
if $A\cap[1,n+1] \in \mathcal S_k$,
then
$A\cap[1,n] \subset A\cap[1,n+1]$ implies that 
$A\cap[1,n] \in \mathcal S_k$.
Thus,
$E_{n+1} \subset E_n$.

We claim that
$$\bigcap_{n=1}^{\infty}E_n\ =\ \mathcal{S}_k.$$
One inclusion is immediate: if $A\in\mathcal S_k$, then by heredity, $A\cap[1,n] \in \mathcal S_k$
for every $n$, so $A\in E_n$ for every $n$. Conversely, suppose that $A\in E_n$ for every $n$. Then
$A\cap[1,n] \in \mathcal S_k$  for every $n$. 
Moreover, $A\cap[1,n] \rightarrow A$
in $\{0,1\}^{\mathbb N}$. By Property (P3), we have $A\in\mathcal S_k$.
This proves the claim.

Since every member of $\mathcal S_k$ is finite, while a random subset
$A$ obtained by independent fair coin flips is almost surely infinite, we have
$$\mathbb P\left(\bigcap_{n=1}^{\infty}E_n\right)\ =\ 0.$$
Since the events $(E_n)_{n=1}^\infty$ are decreasing, continuity of probability from above
gives $\mathbb P(E_n) \rightarrow 0$.
Therefore, 
$$\mathbb P(E_n)\ =\ \frac{|\mathcal S_k\cap\mathcal P([1,n])|}{2^n}\ \ge\ \frac{a_{k,n}}{2^n}$$
implies that 
$$\lim_{n\to\infty}\frac{a_{k,n}}{2^n}\ =\ 0.$$
\end{proof}

\begin{cor}\label{c5}
Let $j\in\mathbb{N}$ and $k\ge 2$. Define 
$$b_{j, k,n}\ =\ 
\left|\left\{F\in j\mathcal S_k\, :\,\max F=n\right\}\right|,\mbox{ for each } n\ge 1.$$
Then 
$$\lim_{n\rightarrow\infty} b_{j,k,n}^{1/n} \ =\ 2\quad\mbox{ and }\quad\lim_{n\to\infty}\frac{b_{j,k,n}}{2^n}\ =\ 0.$$
\end{cor}

\begin{proof}
Since $\mathcal S_k\subset j\mathcal S_k$,
\eqref{e22} gives
$$2^{n-\lceil \log_2 (n+1)\rceil}\ \le\ a_{k,n}\ \le\  b_{j,k,n}\ \le\ 2^{n}.$$
Therefore, 
$$\lim_{n\rightarrow\infty} b^{1/n}_{j,k,n}\ =\ 2.$$

For the second limit, the family $j\mathcal S_k$ is hereditary and consists
only of finite subsets of $\mathbb N$. Furthermore, since $\mathcal S_k$ is compact and the union map from $\mathcal S_k^j$ into $\{0,1\}^{\mathbb N}$ is continuous, $j\mathcal S_k$ is compact.
Thus, the same argument used in the
proof of Proposition \ref{pbound} gives
$$\lim_{n\to\infty}\frac{b_{j,k,n}}{2^n}\ =\ 0.$$
\end{proof}

\begin{proof}[Proof of Corollaries \ref{mc2} and \ref{mc3}]
Corollary \ref{mc2} follows immediately from Theorem \ref{m2} and Proposition \ref{pbound}, while Corollary 1.7 follows from Theorem \ref{m2} and Corollary \ref{c5}.
\end{proof}

%%%%%%%%%%%%%%%%%%%%%%%%%%%%%%%%%%%%%%%%%%%%%%%%%%%%%%%%%%%%%%%%%%%%%%%%%%%%%%%%%%%%%%%%%%%%%%%%%%%%%%%%%%%%%%%%%%%%%%%%
%%%%%%%%%%%%%%%%%%%%%%%%%%%%%%%%%%%%%%%%%%%%%%%%%%%%%%%%%%%%%%%%%%%%%%%%%%%%%%%%%%%%%%%%%%%%%%%%%%%%%%%%%%%%%%%%%%%%%%%%
%%%%%%%%%%%%%%%%%%%%%%%%%%%%%%%%%%%%%%%%%%%%%%%%%%%%%%%%%%%%%%%%%%%%%%%%%%%%%%%%%%%%%%%%%%%%%%%%%%%%%%%%%%%%%%%%%%%%%%%%%

\subsection*{Funding} Not applicable.

\subsection*{Competing interests} The authors declare that they have no competing interests.

\end{document}